\documentclass{ifacconf}

\newtheorem{innercustomdef}{Definition}
\newenvironment{customdef}[1]
  {\renewcommand\theinnercustomdef{#1}\innercustomdef}
  {\endinnercustomdef}  

\newtheorem{innercustomass}{Assumption}
\newenvironment{customass}[1]
  {\renewcommand\theinnercustomass{#1}\innercustomass}
  {\endinnercustomass}

\usepackage{graphics} 
\usepackage{epsfig} 
\usepackage{array}
\usepackage{cite}
\usepackage{algorithm}
\usepackage{acronym}
\usepackage{natbib} 

\usepackage{multicol}
\usepackage{caption}
\usepackage{lipsum}

\acrodef{cav}[\textsc{CAV}]{Connected and Automated Vehicle}
\acrodef{mcav}[\textsc{MCAV}]{Micro Connected and Automated Vehicle}
\acrodef{udssc}[\textsc{UDSSC}]{University of Delaware Scaled Smart City}
\acrodef{ros}[\textsc{ROS}]{Robot Operating System}
\acrodef{fifo}[\textsc{FIFO}]{First-in-First-Out}
\acrodef{soc}[\textsc{SOC}]{State-of-Charge}

\usepackage{amsmath} 
\usepackage{amssymb}  

\begin{document}

\begin{frontmatter}

\title{A Scaled Smart City for Experimental Validation of Connected and Automated Vehicles\thanksref{footnoteinfo}
}

\thanks[footnoteinfo]{This research was supported by the University of Delaware.}     

\author[First]{Adam Stager} 
\author[Second]{Luke Bhan}
\author[First]{Andreas Malikopoulos} 
\author[First]{Liuhui Zhao} 

\address[First]{University of Delaware, Newark, DE 19716 USA (e-mails: astager@udel.edu; andreas@udel.edu; lhzhao@udel.edu).}
\address[Second]{Avon Grove High School, West Grove, PA 19390 USA (e-mail: lukebhan@udel.edu)}


\begin{abstract}
The common thread that characterizes energy-efficient mobility systems for smart cities is their interconnectivity which enables the exchange of massive amounts of data. This, in turn, provides the opportunity to develop a decentralized framework to process this information and deliver real-time control actions that optimize energy consumption and other associated benefits. To seize these opportunities, this paper describes the development of a scaled smart city providing a glimpse that bridges the gap between simulation and full scale implementation of energy-efficient mobility systems. Using this testbed, we can quickly, safely, and affordably experimentally validate control concepts aimed at enhancing our understanding of the implications of next generation mobility systems.
\end{abstract}


\end{frontmatter}


\section{Introduction}
In a rapidly urbanizing world, we need to make fundamental transformations in how we use and access transportation.  Energy-efficient mobility systems such as \acp{cav} along with shared mobility and electric vehicles provide the most intriguing and promising opportunities for enabling users to better monitor transportation network conditions and make better operating decisions to reduce energy consumption, greenhouse gas emissions, travel delays and improve safety. As we move to increasingly complex transportation systems new control approaches are needed to optimize the system behavior resulting from the interactions between vehicles navigating different traffic scenarios. 

Given this new environment, the overarching goal of this paper is to (1) report on the development of the \ac{udssc} testbed that includes 35 robotic cars to replicate real-world traffic scenarios in a small and controlled environment, and (2) use this testbed to demonstrate CAV coordination at merging roadways. \ac{udssc} can serve as a testbed to explore the acquisition and processing of vehicle-to-vehicle and vehicle-to-infrastructure communication. It can also help us prove control concepts on coordinating \acp{cav} in specific transportation scenarios, e.g., intersections, merging roadways, roundabouts, speed reduction zones, etc. These scenarios along with the driver responses to various disturbances are the primary sources of bottlenecks that contribute to traffic congestion; see \cite{Margiotta2011,Malikopoulos2013}. 
In 2015, congestion caused people in urban areas in the US to spend 6.9 billion additional hours on the road and to purchase an extra 3.1 billion gallons of fuel, resulting in a total cost estimated at \$160 billion; see \cite{Schrank2015}.  

\acp{cav} can provide shorter gaps between vehicles and faster responses while improving highway capacity. Several research efforts have been reported in the
literature proposing either \emph{centralized} or \emph{decentralized} approaches
for coordinating \acp{cav} in specific traffic scenarios. The overarching goal of such efforts is to yield a smooth traffic flow avoiding stop-and-go driving. Numerous approaches have been reported in the literature on coordinating \acp{cav} in different transportation scenarios with the intention of improving traffic flow. \cite{Kachroo1997} proposed a longitudinal and lateral controller to guide the vehicle until the merging maneuver is completed. Other efforts have focused on developing a hybrid control aimed at keeping a safe headway between vehicles in the merging process, see \cite{Kachroo1997, Antoniotti1997}; or developing three levels of assistance for the merging process to select a safe space for the vehicle to merge; see \cite{Ran1999}. Some authors have explored virtual vehicle platooning, where a controller identifies and interchanges appropriate information between the vehicles involved in the merging maneuver while each vehicle assumes its own control actions to satisfy the assigned time and reference speed; see \cite{Lu2000}.

\cite{VanMiddlesworth2008} addressed the problem of traffic coordination for small intersections which commonly handle low traffic loads. \cite{Milanes2010} designed a controller that allows a fully automated vehicle to yield to an incoming vehicle in the conflicting road or to cross, if it is feasible without the risk of 
collision. \cite{Alonso2011} proposed two conflict resolution schemes in which an autonomous vehicle could make a decision about the appropriate crossing schedule and trajectory to follow to avoid collision with other manually driven vehicles on the road.  A survey of the research efforts in this area that have been reported in the literature to date can be found in \cite{Malikopoulos2016a}.

Although previous work has shown promising results emphasizing the potential benefits of coordination between \acp{cav}, validation has been primarily in simulation. In this paper, we demonstrate coordination of scaled \acp{cav} and quantify the benefits in energy usage. The contributions of this paper are: 1) the development of a 1:25 scaled smart city capable of testing decentralized control algorithms on up to 35 \acp{mcav} and 2) the experimental validation of a control framework reported in \cite{Rios-Torres2} for coordination of \acp{cav}.

The remainder of the paper proceeds as follows. In Section II, we introduce the configurations of \ac{udssc}. In Section III, we review a decentralized control framework for coordination of \acp{cav} in merging roadways. Experimental results in Section IV illustrate the effectiveness of the proposed solution in the scaled smart city environment. We draw concluding remarks in Section V.

\section{University of Delaware Scaled Smart City (UDSSC)} \label{scaledsmartcity}

\begin{figure}
  \centering
    \includegraphics[width=3.2 in]{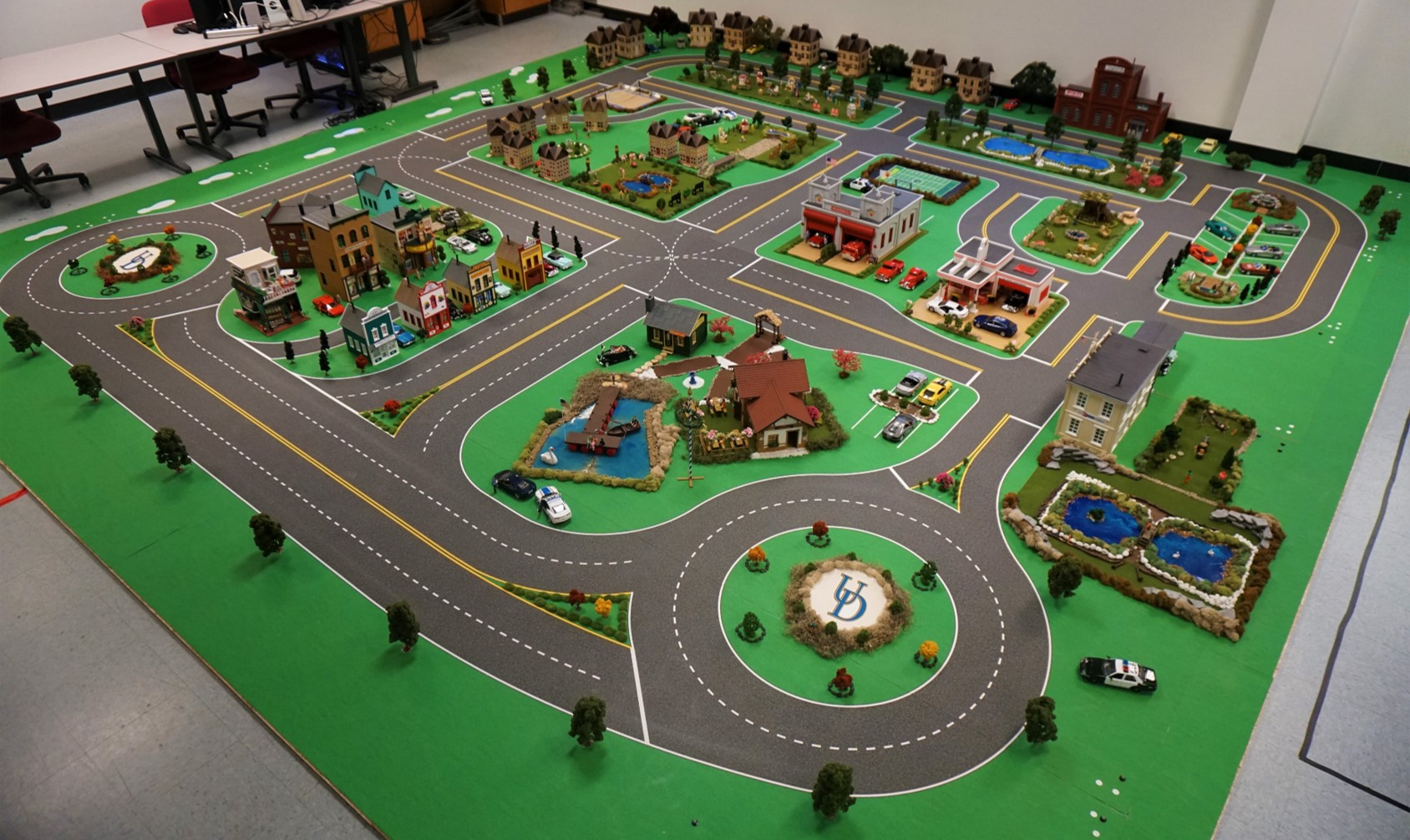} 
      \caption{Birdseye view of the \ac{udssc}.}
      \label{udssc}
\end{figure}

\ac{udssc} is a testbed that can replicate real-world traffic scenarios in a small, controlled environment and help formulate the appropriate features of a ``smart" city. It can be used as an effective way to visualize the concepts developed using \acp{cav} and their related implications in energy usage. \ac{udssc} is a fully integrated smart city (Fig. \ref{udssc}) 
incorporating realistic environmental cues, scaled \acp{mcav}, and state-of-the-art, high-end computers supporting standard software for system analysis and optimization for simulating different control algorithms and distributing control inputs to as many as 35 \acp{mcav}.

\subsection{Physical Design I: Map}
The \ac{udssc} spans over 400 square feet and includes one lane intersections, two lane intersections, roundabouts, and a highway with entrance and exit ramps. Using SolidWorks to accurately maintain 1:25 scale, a fully dimensioned blueprint 
was designed forming a cohesive roadway representative of real world road scenarios. 
Using a HP DesignJet z5200 printer the map is printed with realistic texturing on twelve 44"x120" sheets of wear resistant HP Artist Matte Canvas that can be easily replaced to either reconfigure or repair sections of the city independently. 
Eight Vicon Vantage V16 cameras are used to localize the map within a global coordinate system. 


\subsection{Physical Design II: Cars}
\begin{figure}
  \centering
    \includegraphics[width=3.3 in]{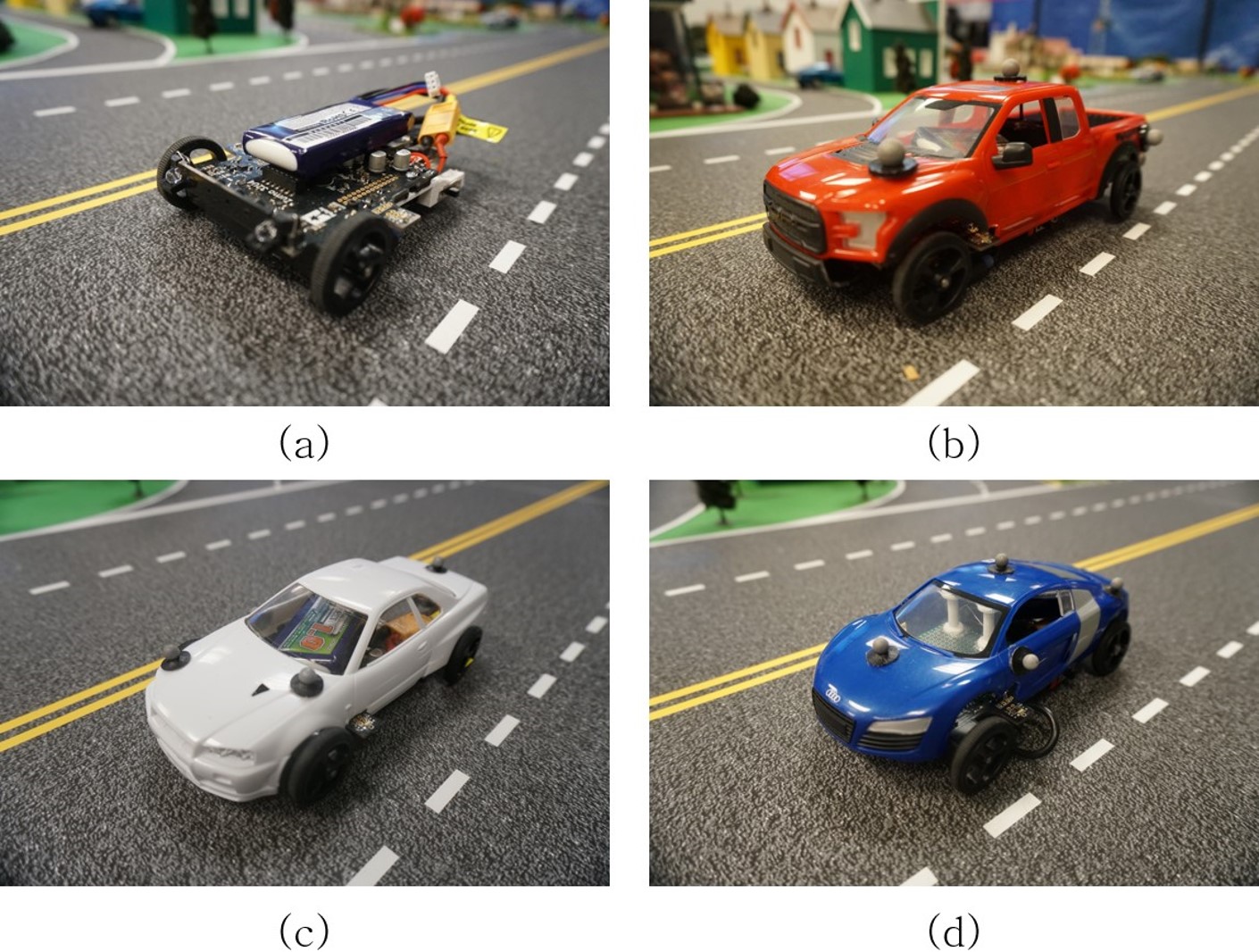} 
      \caption{The \ac{udssc} incorporates 35 automated vehicles, starting with a \ac{mcav} platform (a) and incorporating different types of scaled vehicle shells such as a Ford F150 (b), Nissan Skyline (c), and Audi R8 (d).}
      \label{mcavs}
\end{figure}

Scaled \acp{mcav} have been designed using easily assembled off-the-shelf components coupled with several 3D printed parts (Fig. \ref{mcavs}). At the core of each platform is a 75.81:1 geared, differentially driven Pololu Zumo, offering dual H-bridge motor drivers, $n_e = 12$ counts per revolution (CPR) encoders and an on-board Atmega 32U4 micro-controller. Each Zumo contains an embedded set of sensors including an IMU, line-following, and infrared proximity sensors which can provide feedback to each \ac{mcav}. Rubberized wheels with radius $r=1.6$ $cm$ are mounted directly to each gear motor output shaft and separated by $d=9$ $cm$ to roughly mimic the 1:25 scale width of full sized cars/trucks. The Zumo is connected to an on-board Raspberry Pi 3 with 1.2 GHz quad-core ARM Cortex A53 micro-processor and WiFi used for communication. The \ac{mcav} platform (not including its car-shaped shell) measures $13$ $cm$ x $10.5$ $cm$ x $4.5$ $cm$ $(l/w/h)$. 
A power regulator manages the voltage requirements of the Raspberry Pi 3, supplying a regulated 5VDC from a 7.4VDC, 1000mAh Li-Po battery. Fully charged \acp{mcav} are capable of running approximately 90 minutes before being recharged.

\subsection{Control System Architecture}

\ac{udssc} has a multi-level control architecture 
with high level commands originating from a centralized PC called the ``Main Frame" (Processor: Intel Core i7-6950X CPU @ 3.00 GHz x 20, Memory: 125.8 GiB) then enforced by a low-level controller on-board each \ac{mcav}. Asynchronous communication is enabled by WiFi connectivity using the Node.js non-blocking I/O protocol. Electron uses the socketIO library supporting multi-threading for multi-vehicle communication. A web browser combines JavaScript, HTML, and CSS for a user friendly interface into the \ac{ros} architecture. Generally control of each \ac{mcav} can be broken into merging, lane and reference tracking controllers, the latter two of which enforce realistic road behaviors (i.e. staying in the center of the road and respecting speed limits.)


\subsubsection{Lane tracking:}
The roads of \ac{udssc} are encoded as sequences of tangent arcs and straight line segments each with an associated potential field. 
Potential field methods 
provide important computational expedience due to their analytical representation. Given a looping sequence of road segments an \ac{mcav} starting on the first road segment will follow along the road until its battery is depleted.
Although an \ac{mcav} can use lane tracking directly to follow along a road, while in the control region merging requires a carefully maintained velocity profile with respect to the lane center. Offsets from the centerline of the lane using lane tracking results in noisy forward velocity measurements, however a virtual robot can track the center of each lane exactly. For situations where precise velocity profiles are required, instead of controlling each \ac{mcav} using the potential field directly a virtual robot is simulated within the vector field and used as a reference point tracked by the real robot.

\subsubsection{Reference tracking:}
As long as a reference point tracks a sequence of road segments using lane tracking, an \ac{mcav} can be controlled by reference tracking. 
A state tracking method as described in \cite{Giuseppe2002} is used
, approximately linearizing the error dynamics of the \ac{mcav}'s local frame with respect to a reference trajectory. 

\subsubsection{Merging scenario:}
There have been several approaches for automated vehicle merging as reported in \cite{Malikopoulos2016a}. In this paper, we consider the control approach described in Section \ref{sec:control} with the Main Frame tracking the positions of each \ac{mcav} in order to determine when vehicles enter the control zones. 
Once inside the control zone a virtual robot tracks the desired velocity profile exactly and reference tracking is used to mimic this behavior by the associated \ac{mcav}. 

\subsubsection{Low level control:}
In place of car-like models which would require a more complex mechanical design, Zumo based \acp{mcav} are differentially driven. Practically the control input is expressed as $u = [v,\omega]$ with $v$ and $\omega$ representing the forward and angular velocities respectively. On-board encoders enable low-level control from wirelessly transmitted inputs, converted into independent wheel velocities by two relationships: $v = R(\dot{\phi}_R+\dot{\phi}_L)/2$ and $\omega = R(\dot{\phi}_R-\dot{\phi}_L)/2$ with
$\dot{\phi}_R,\dot{\phi}_L$ 
as the right and left wheels angular velocity respectively and $R$ as wheel radius. High frequency control results in noisy measurements due to low-resolution encoders. The Atmega 32U4 measures encoder pulses at a frequency of 2 kHz then smooths the velocity estimate by averaging a 25 measurement queue. 
A proportional controller adjusts the PWM duty cycle depending on the error between measured and desired velocity, saturating at 0.7m/s $\pm$ 0.1m/s depending on transmission friction and slight variations between \acp{mcav}.

\begin{figure}
  \centering
    \includegraphics[width=3.3
    in]{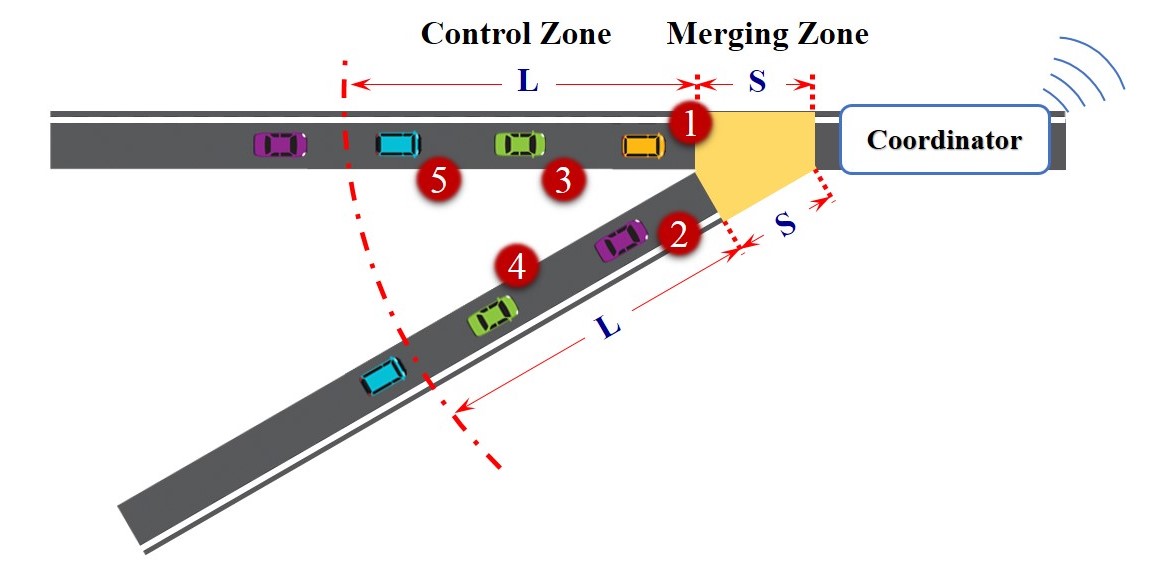} 
      \caption{Merging roads with \acp{cav}.}
      \label{fig:merging}
\end{figure}

\section{Coordination of Connected and Automated Vehicles} \label{sec:control}
\subsection{Modeling Framework}
We consider a merging roadway (Fig. \ref{fig:merging}) consisting of main and secondary roads. The region where lateral collision between vehicles can occur is called \textit{merging zone} and has a length of $S$. On each road, there is a \textit{control zone} inside of which all vehicles can communicate with each other and with a coordinator. Note that the coordinator is not involved in any decision for any \ac{cav} and only enables communication of appropriate information among \acp{cav}. The distance from the entry of the control zone to the entry of the merging zone is $L$. The value of $L$ depends on the communication range between \acp{cav} and the coordinator, and $S<L$ is the physical length of a merging zone. 

Let $N(t)\in\mathbb{N}$ be the number of \acp{cav} inside the control zone at time
$t\in\mathbb{R}^{+}$ and $\mathcal{N}(t)=\{1,\ldots,N(t)\}$ be a queue designating the order in which these vehicles enter the merging zone. Thus,
letting $t_{i}^{m}$ be the assigned time for vehicle $i$ to enter the merging zone, we
require that
\begin{equation}
t_{i}^{m}\geq t_{i-1}^{m},~\forall i\in\mathcal{N}(t),~i>1. \label{eq:fifo}%
\end{equation}
There are a number of ways to satisfy \eqref{eq:fifo}. For example, we may
impose a strict \ac{fifo} queuing structure, where each vehicle must
enter the merging zone in the same order it entered the control zone. 
We investigate a specific scheme for determining
$t_{i}^{m}$ (upon arrival of \ac{cav} $i$) based on our problem formulation, without affecting $t_{1}^{m},\ldots,t_{i-1}^{m}$, but emphasize that our analysis is not restricted by the policy designating the order of the vehicles
within the queue $\mathcal{N}(t)$.

We adopt the optimization framework proposed in \cite{Rios-Torres2} for coordinating the merging of \acp{cav}. 
The dynamics of each
vehicle $i\in\mathcal{N}(t)$ are represented by a double integrator,
\begin{equation}
\dot{p}_{i} =v_{i}(t), ~\dot{v}_{i} =u_{i}(t), \label{eq:model}
\end{equation}
where $t\in\mathbb{R}^{+}$ is the time, $p_{i}(t)\in\mathcal{P}_{i}$, $v_{i}(t)\in\mathcal{V}_{i}$,
and $u_{i}(t)\in\mathcal{U}_{i}$ denotes position, speed and acceleration/deceleration
(control input) of each vehicle $i\in\mathcal{N}(t)$ inside the control zone. Let
$\left[\begin{array}{cc}
p_{i}(t) & v_{i}(t)\end{array}\right]^{T}$ denote the state of each vehicle $i$, with initial value $\left[\begin{array}{cc}
0 & v_{i}^{0}\end{array}\right]^{T}$.  The state space $\mathcal{P}%
_{i}\times\mathcal{V}_{i}$ is closed with respect to the induced topology,
thus, it is compact.

For any initial state  $\left[\begin{array}{cc}
p_{i}(t_{i}^{0}) & v_{i}(t_{i}^{0})\end{array}\right]^{T}$, where $t_{i}^{0}$ is the time that the vehicle $i$ enters the control zone,
and every admissible control $u(t)$, the double integrator
has a unique solution on some interval $[t_{i}^{0},t_{i}^{m}]$,
where $t_{i}^{m}$ is the time that vehicle $i\in\mathcal{N}(t)$
enters the merging zone. In our framework we impose the following state and control constraints:
\begin{equation}%
\begin{split}
u_{i,min} &  \leqslant u_{i}(t)\leqslant u_{i,max},\quad\text{and}\\
0 &  \leqslant v_{min}\leqslant v_{i}(t)\leqslant v_{max},\quad\forall
t\in\lbrack t_{i}^{0},t_{i}^{m}],
\end{split}
\label{speed_accel constraints}%
\end{equation}
where $u_{i,min}$, $u_{i,max}$ are the minimum and maximum
control inputs (maximum deceleration/acceleration) for each vehicle $i\in\mathcal{N}(t)$, and $v_{min}$, $v_{max}$
are the minimum and maximum speed limits respectively. For simplicity, in the
rest of the paper we consider no vehicle diversity, and thus, we set $u_{i,min}=u_{min}$
and $u_{i,max}=u_{max}$.

For absence of any rear-end collision of two consecutive vehicles
traveling on the same lane, the position of the preceding vehicle should be
greater than or equal to the position of the following vehicle plus a safe distance $\delta(v_{ave}(t))<S$, which is a function of the average speed of the vehicles inside the control zone. Thus, we impose the following rear-end safety
constraint
\begin{equation}
s_{i}(t)=p_{k}(t)-p_{i}(t)\geqslant\delta(v_{ave}(t)),~\forall t\in\lbrack t_{i}^{0}%
,t_{i}^{m}],\label{eq:rearend}%
\end{equation}
where $k$ denotes the vehicle that is physically
located ahead of $i$ in the same lane, and $v_{ave}(t)$ is the average speed of the vehicles inside the control zone at time $t$.


\begin{customdef}{1}
Each CAV $i\in\mathcal{N}(t)$ belongs to at least one of the following two
subsets of $\mathcal{N}(t)$ depending on its physical location inside the control zone:
1) $\mathcal{L}_{i}(t)$ contains all \acp{cav} traveling on
the same road and lane as vehicle $i$ and 2) $\mathcal{C}_{i}(t)$ contains all \acp{cav} traveling on a different
road from $i$ and can cause collision at the merging zone. \label{def:1}
\end{customdef}

\begin{customdef}{2}
For each vehicle $i\in\mathcal{N}(t)$, we define the set $\Gamma_{i}$ that
includes only the positions along the lane where a lateral collision is
possible, namely
\begin{equation}
\Gamma_{i}\triangleq\Big\{t~|~t\in\lbrack t_{i}^{m},t_{i}^{f}]\Big\},
\label{eq:def1}
\end{equation}
with $t_{i}^{f}$ as the time vehicle $i$ exits the merging zone.
\end{customdef}

Consequently, to avoid a lateral collision for any two vehicles $i,j\in
\mathcal{N}(t)$ on different roads, the following constraint should hold
\begin{equation}
\Gamma_{i}\cap\Gamma_{j}=\varnothing,\text{ \ \ \ }\forall t\in\lbrack
t_{i}^{m},t_{i}^{f}]\text{, \ }j\in\mathcal{C}_{i}(t). \label{eq:lateral}%
\end{equation}

The above constraint implies that only one vehicle at a time can be inside the
merging zone. If the length of the merging is long, then this constraint may not be
realistic since it results in dissipating space and capacity of the road.
However, the constraint is not restrictive in the problem formulation and it
can be modified appropriately.

In the modeling framework described above, we impose the following assumptions: 


\begin{customass}{1} \label{ass:srz} The vehicles cruise inside the merging zone with an imposed speed limit, $v_{srz}$.
This implies that for each vehicle $i$
\begin{equation}
t_{i}^{f} = t_{i}^{m}+\frac{S}{v_{srz}}.\label{eq:time}%
\end{equation}
\end{customass}

This assumption is intended to enhance safety awareness, but it could be modified appropriately, if necessary.


\subsection{Communication Structure of Connected and Automated Vehicles}
We consider the problem of deriving the optimal control input (acceleration/deceleration)
of each \ac{cav} inside the control zone (Fig. \ref{fig:merging}), under the hard safety
constraints to avoid rear-end and lateral collision. By controlling the speed of the vehicles, the speed of the queue built-up at the merging zone decreases, and thus the congestion recovery time is also reduced. The latter results in maximizing the throughput in the merging zone.

When a \ac{cav} $i$ enters the control zone, it can communicate with the other \acp{cav} in the control zone and with the coordinator. Note that the coordinator is not involved in any decision for any \ac{cav} and only enables communication of appropriate information among \acp{cav}. The coordinator handles the information between the vehicles as follows. When a \ac{cav} reaches the control zone  at some instant $t$, the coordinator assigns a \textit{unique identity} to
each vehicle $i\in\mathcal{N}(t),$ which is a pair $(i,j)$, where $i=N(t)+1$ is an integer representing the location of the vehicle in a \ac{fifo} queue $\mathcal{N}(t)$ and $j\in\{1,2\}$ is an integer 
based on a one-to-one mapping from 
$\mathcal{L}_i(t)$ and $\mathcal{C}_i(t)$ onto
$\{1,2\}$. If the vehicles enter the control zone at the same time, then the coordinator selects randomly their position in the queue. 

\begin{customdef}{3} For each \ac{cav} $i$ entering the control zone, we define the \textit{information set}
$Y_{i}(t),$ which includes all information that each vehicle shares, as
\begin{gather}
Y_{i}(t) \triangleq\Big\{p_{i}(t),v_{i}(t),\mathcal{Q},t_{i}%
^{m}\Big\},\forall t \in\lbrack t_{i}^{0},t_{i}^{m}],
\end{gather}
\label{infoset}where $p_{i}(t),v_{i}(t)$ are the position and speed of \ac{cav} $i$
inside the control zone, $\mathcal{Q}\in\{\mathcal{L}_i(t),\mathcal{C}_i(t)\}$ is the subset assigned to \ac{cav} $i$ by the coordinator, and $t_{i}^{m}$, is the target time for \ac{cav} $i$ to
enter the merging zone.
\end{customdef}

The time $t_{i}^{m}$ that vehicle $i$ will be entering the merging zone is restricted by imposing rear-end and lateral collision constraints. Therefore, to ensure \eqref{eq:rearend} and \eqref{eq:lateral} are satisfied at $t_{i}^{m}$ we impose the following conditions which depend on the subset vehicle $i-1$ belongs to.
If $i-1\in\mathcal{L}_{i}(t),$ 
\begin{equation}
t_{i}^{m}=\max\Bigg\{ \min\Big\{t_{i-1}^{m}+\frac{\delta(v_{ave}(t))}{v_{srz}},\frac{L}{v_{min}}\Big\}, \frac{L}{v_{i}(t_{i}^{0})}, \frac{L}{v_{max}}\Bigg\},\label{eq:condition1a}%
\end{equation}
and if $i-1\in\mathcal{C}_{i}(t),$ 
\begin{equation}
t_{i}^{m}=\max\Bigg\{ \min\Big\{t_{i-1}^{m}+\frac{S}{v_{srz}},\frac{L}{v_{min}}\Big\}, \frac{L}{v_{i}(t_{i}^{0})}, \frac{L}{v_{max}}\Bigg\},\label{eq:condition1b}%
\end{equation}

where $v_{srz}$, is the imposed speed inside the merging zone (Assumption \ref{ass:srz}), and $v_{i}(t_{i}^{0})$ is the initial speed of vehicle $i$ when it enters the control zone at $t_{i}^{0}$. The conditions \eqref{eq:condition1a} and \eqref{eq:condition1b} ensure the time $t_{i}^{m}$ each vehicle $i$ will be entering the merging zone is feasible and can be attained based on the imposed speed limits inside the control zone. In addition, for low traffic flow where vehicle $i-1$ and $i$ might be located far away from each other, there is no compelling reason for vehicle $i$ to accelerate within the control zone to maintain a distance $\delta(v_{ave}(t))$ from vehicle $i-1$. Therefore, in such cases vehicle $i$ can keep cruising within the control zone with the initial speed $v_{i}(t_{i}^{0})$ that entered the control zone at $t_{i}^{0}.$ 

The recursion is initialized when the first vehicle enters the control zone and is assigned $i=1$. In this case, $t_{1}^{m}$ can be externally assigned as the desired exit time of this vehicle whose behavior is
unconstrained. Thus the time $t_{1}^{m}$ is 
fixed and available through $Y_{1}(t)$. The
second vehicle accesses $Y_{1}(t)$ to compute time $t_{2}^{m}$. The third vehicle accesses $Y_{2}(t)$ and the communication process continues in the same fashion until vehicle $N(t)$ in the queue accesses $Y_{N(t) -1}(t)$.

\subsection{Optimal Control Problem Formulation for \acp{cav}}
Since the coordinator is not involved in any decision on the vehicle
coordination we can formulate $N(t)$ sequential
decentralized control problems that may be solved on-line:
\begin{gather}
\min_{u_{i}}\frac{1}{2}\int_{t_{i}^{0}}^{t_{i}^{m}}u_{i}^{2}%
(t)~dt,\label{eq:decentral}\\
\text{subject to}:\eqref{eq:model}~\text{and}~ \eqref{speed_accel constraints},\nonumber
\end{gather} 
with initial and final conditions: $p_{i}(t_{i}^{0})=0$, $p_{i}(t_{i}^{m})=L,$ $t_{i}^{0}$, $v_{i}(t_{i}^{0}),$ $t_{i}^{m},$ and $v_{i}(t_{i}^{m})=v_{srz}$. In \eqref{eq:decentral}, rear end
 \eqref{eq:rearend} and lateral \eqref{eq:lateral} collision safety constraints are omitted. As mentioned earlier,  \eqref{eq:lateral} implicitly handled by
the selection of $t_{i}^{m}$ in \eqref{eq:condition1b}. Eq. \eqref{eq:rearend} is omitted because it has been shown \cite{Malikopoulos2017} that the solution of
\eqref{eq:decentral} guarantees this constraint holds throughout $[t_{i}^{0},t_{i}^{f}]$. Thus, \eqref{eq:decentral} is a simpler problem to solve on-line. 

The analytical solution of \eqref{eq:decentral} without state and control constraints was presented in
 \cite{Rios-Torres2015, Rios-Torres2, Ntousakis:2016aa} for real time coordination of \acp{cav} at highway on-ramps and \cite{ZhangMalikopoulosCassandras2016} at two adjacent intersections. 

\section{Experimental Results}
\begin{figure}[!ht]
  \centering
    \includegraphics[width=3.3 in]{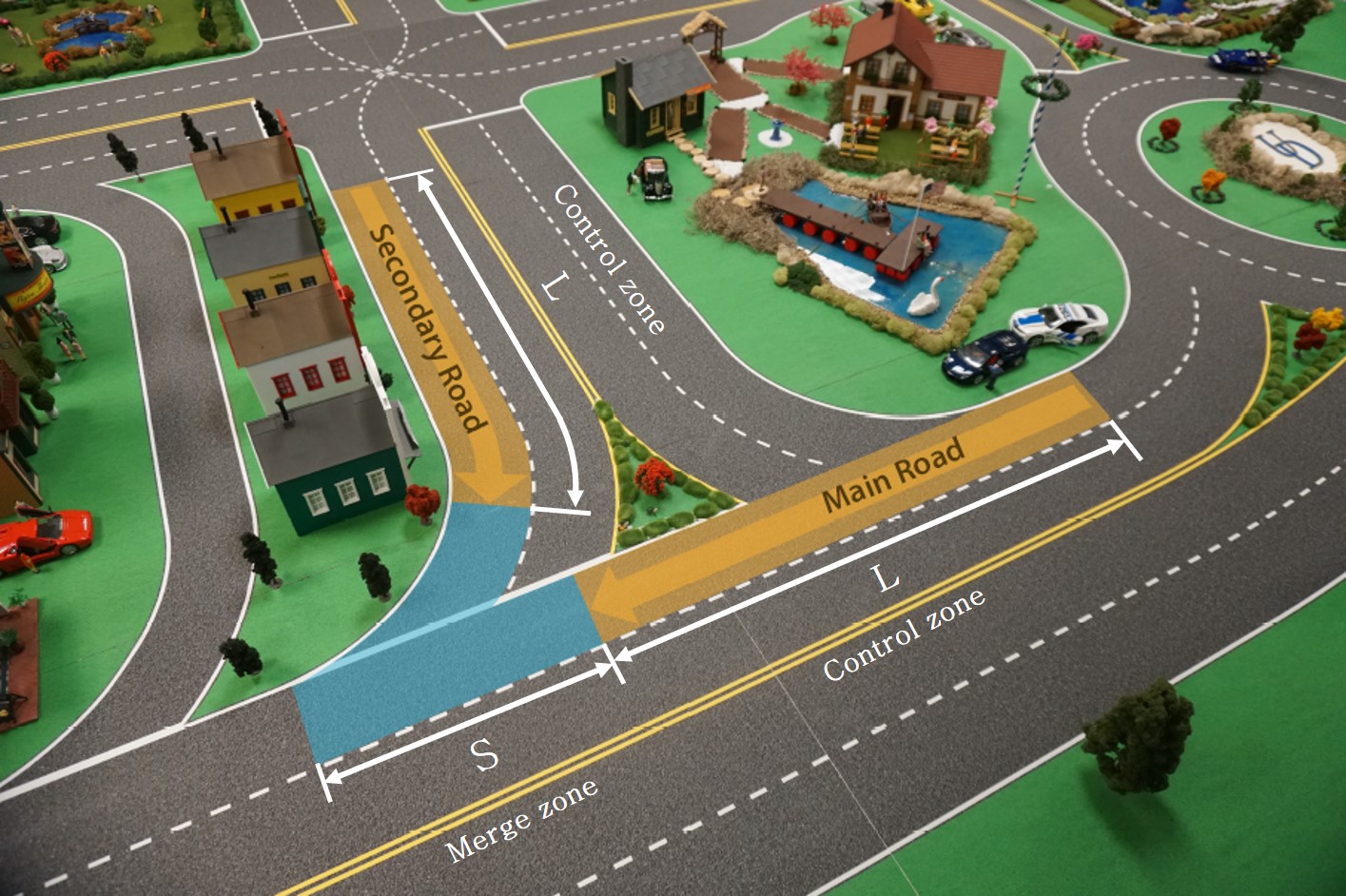} 
      \caption{Aerial view of real control and merging regions.}
 \label{fig:experiment_setup}
\end{figure}

To evaluate the effectiveness and efficiency of the proposed approach, a total number of ten \acp{mcav} are set up in a merging scenario  [video available in \cite{Stager2017}]. Five \acp{mcav} cruise on the main road in \ac{udssc}, while the other five \acp{mcav} cruise on the secondary road with the intention to merge into the main road (Fig.~ \ref{fig:experiment_setup}).

We consider two scenarios: 1) all \acp{mcav} are controlled by the decentralized control algorithm; 2) all \acp{mcav} behave based on a simple (baseline) car following model, with \acp{mcav} on the secondary road yielding to \acp{mcav} of the main road to avoid lateral collision in the merging zone.


\begin{figure}[!ht]
  \centering
    \includegraphics[width=3.3 in]{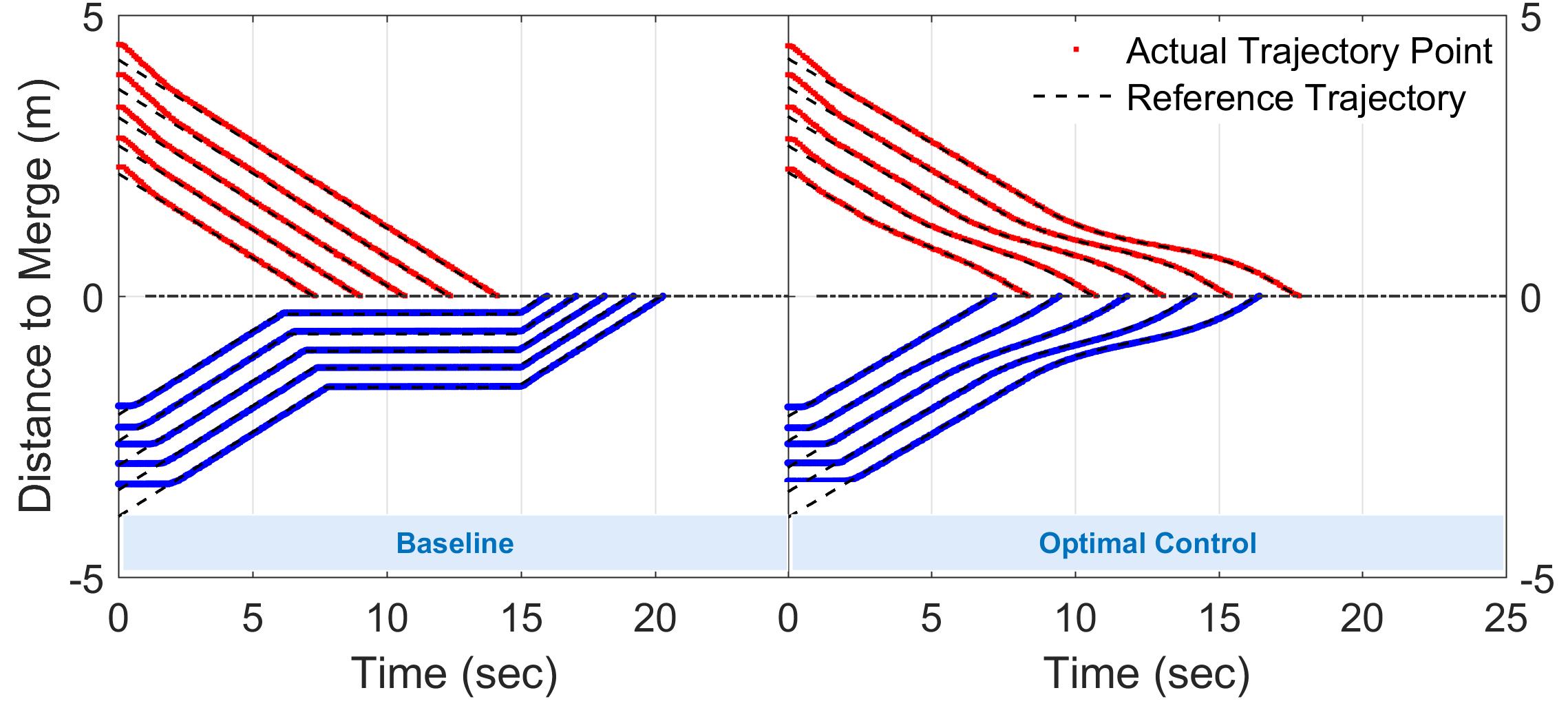} 
      \caption{Comparison of vehicle trajectories.}
 \label{fig:traj}
\end{figure}

\subsubsection{Vehicle Position Trajectory:}
Position trajectories of \acp{mcav} under two scenarios are illustrated in Fig.~\ref{fig:traj}. A dashed line represents the reference trajectory for each vehicle commanded by the control algorithm, while a dense scatter plot represents points measured along the actual trajectory achieved by each \ac{mcav}. To separate the \acp{mcav} on two roads, trajectories are flipped over Y-axis. Thus, in Fig.~\ref{fig:traj}, red dots stand for trajectory points of the \acp{mcav} on the main road, and blue dots stand for the trajectory points of \acp{mcav} on the secondary road. On the right panel of Fig.~\ref{fig:traj}, \acp{mcav} follow the optimal trajectory and merge successfully without stop-and-go driving with only marginal errors. Position trajectories of \acp{mcav} cruising without the optimal control (baseline scenario) are shown in the left panel of Fig.~\ref{fig:traj}. Since the gaps between the mainline cars are not large enough for merging cars to safely merge into the roadway, merging cars need to stop until all leading mainline vehicles traverse the merging zone, resulting in queuing on the merging roadway. For comparison, the merging maneuver for all ten cars is completed in 16.5 sec with the optimal control algorithm, whereas it takes 20.3 sec for the baseline (i.e. an 18.7\% travel time savings is achieved with the optimal control algorithm.)

\subsubsection{Battery State-of-Charge:}
\begin{figure}
  \centering
    \includegraphics[width=3.3 in]{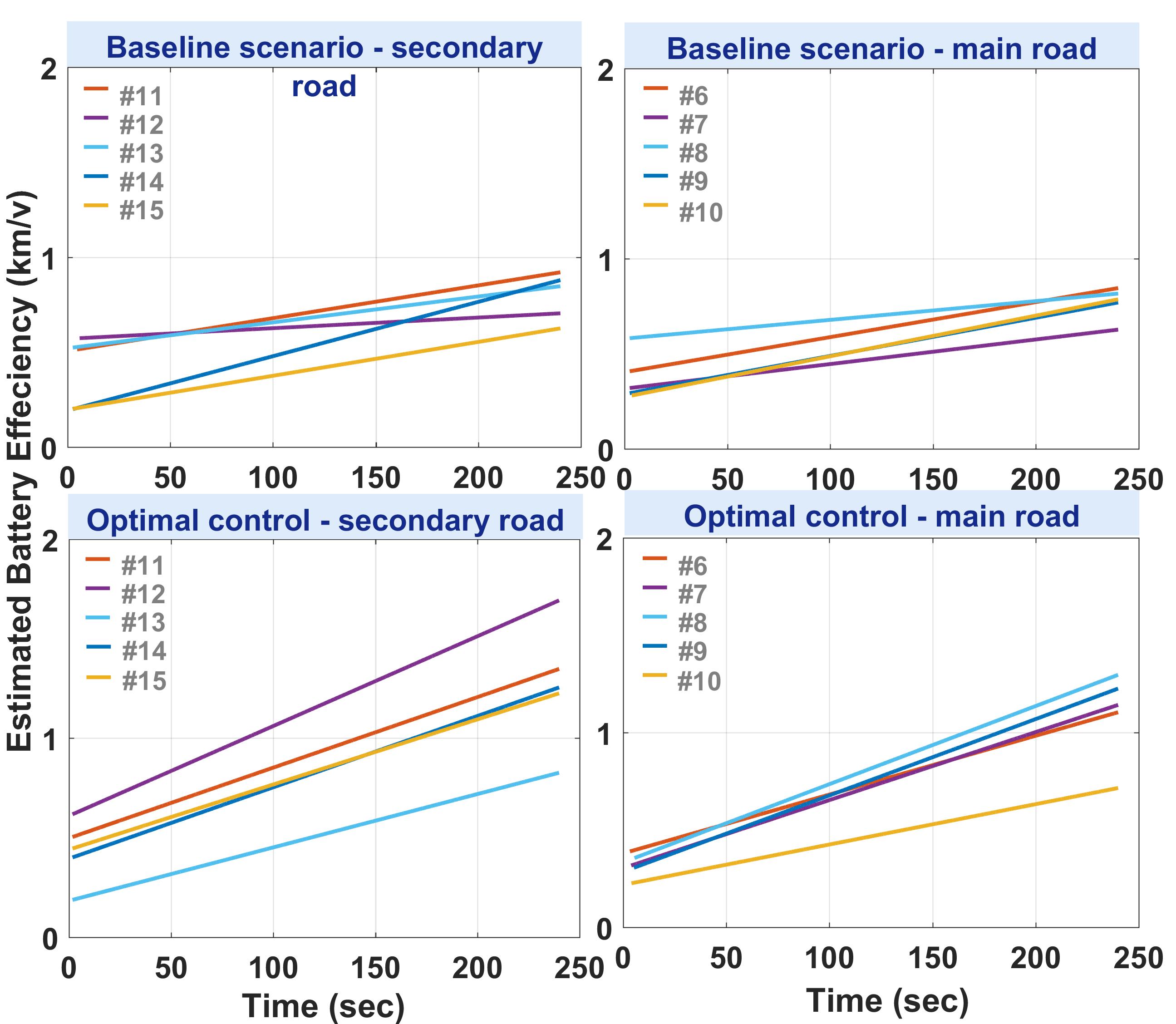} 
      \caption{Battery efficiency over time.}
 \label{fig:veh_battery}
\end{figure}

\begin{figure}
  \centering
    \includegraphics[width=3.3 in]{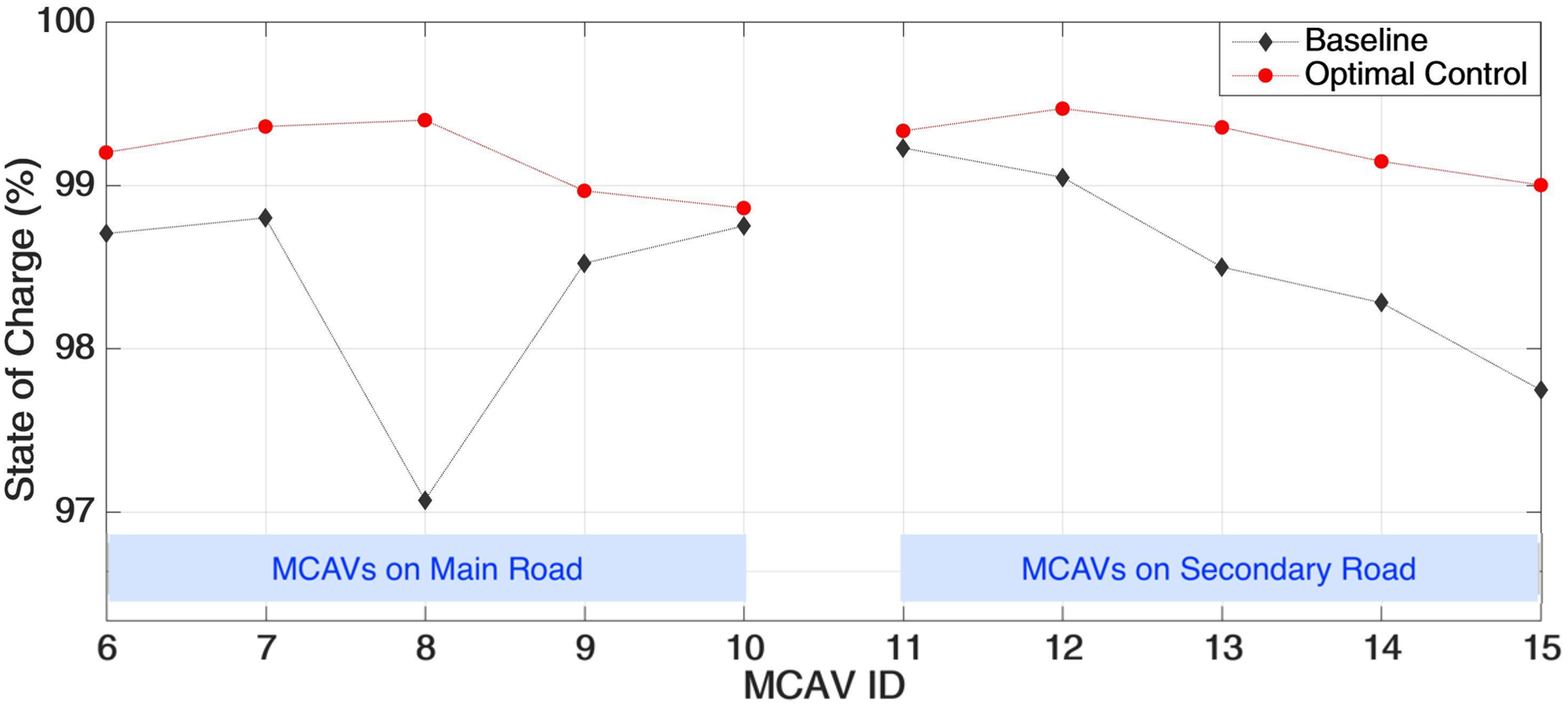} 
      \caption{Final state of charge of the battery for each \ac{mcav}.}
 \label{fig:soc}
\end{figure}

To quantify the benefits of vehicle coordination, we compare the battery \ac{soc} for each \ac{mcav}. Under both scenarios, the \acp{mcav} loop around the merging zone following a predefined trajectory. \ac{soc} is recorded for a 4-minute run. The estimated battery efficiencies for \acp{mcav} under the two scenarios are illustrated in Fig.~\ref{fig:veh_battery}. From the final \ac{soc} of each \ac{mcav} (Fig.~\ref{fig:soc}), it is clear that coordination of \acp{mcav} improves the energy efficiency in the merging scenario due to the elimination of the stop-and-go driving. 

\section{Concluding Remarks}
\ac{udssc} is a small-scale ``smart" city that can replicate real-world traffic scenarios in a small and controlled environment. This testbed can be an effective way to visualize the concepts developed in real world traffic scenarios using \acp{cav} in a quick, safe, and affordable way. The \ac{udssc} helps bridge the gap between theory and practical implementation by providing a means of simultaneously testing as many as 35 \acp{mcav}. We used \ac{udssc} to validate experimentally a control framework reported in \cite{Rios-Torres2} for coordination \acp{cav}. The results demonstrate that coordination of \acp{cav} can improve the battery efficiency due to elimination of the stop-and-go driving. The integration of human driven \acp{mcav} adds a promising direction for future work that may provide insights on the impact of \acp{cav} in real world scenarios. 

\begin{ack}
Caili Li (Univ. of Delaware) for assistance with software, Yue Feng (Univ. of Delaware) for WiFi communication and user interface, Grace Gong (Wilmington Charter High School) for \ac{mcav} design, and Abhinav Ratnagiri (Concord High School) for serial communication. 
\end{ack}
\bibliography{TITS_merging}

\end{document}